\documentclass[12pt,reqno]{amsart}
\usepackage[margin=0.95in]{geometry}
\usepackage{amscd,amsfonts,amsmath,amssymb,amsthm,latexsym,mathrsfs,textcomp,verbatim}
\usepackage{accents}
\usepackage{bm}
\usepackage{bbm}
\usepackage{booktabs}
\usepackage{cite}
\usepackage{color}
\usepackage{constants}
\usepackage{tikz}
\usetikzlibrary{positioning, arrows.meta}
\usepackage{csquotes}
\usepackage{enumerate}
\usepackage{etoolbox}
\usepackage{float}
\usepackage{hypbmsec}
\usepackage[dvipsnames]{xcolor}
\usepackage[breaklinks=true,colorlinks=true,linkcolor=black!40!blue,citecolor=black!40!blue,filecolor=black!40!blue ,urlcolor=black!40!blue, pagebackref]{hyperref}
\hypersetup{linktocpage}

\usepackage[capitalise]{cleveref}
\usepackage{mathtools}
\usepackage{soul}
\usepackage{tikz}
\usepackage{tikz-cd}
\usetikzlibrary{shapes.geometric}
\usetikzlibrary{shapes.misc}
\usetikzlibrary{positioning}
\allowdisplaybreaks[4]
\usepackage{hypbmsec}
\newtheorem{theorem}{Theorem}[section]

\newtheorem{conjecture}[theorem]{Conjecture}

\newtheorem{remark}[theorem]{Remark}
\newtheorem{definition}[theorem]{Definition}

\theoremstyle{definition}

\numberwithin{equation}{section}

\newcommand{\R}{\mathbb{R}}
\newcommand{\Q}{\mathbb{Q}}

\renewcommand{\Re}{\mathrm{Re}}
\renewcommand{\leq}{\leqslant}
\renewcommand{\geq}{\geqslant}
\usepackage{orcidlink}

\patchcmd{\section}{\scshape}{\bfseries}{}{}
\makeatletter
\renewcommand{\@secnumfont}{\bfseries}
\makeatother

\makeatletter\newcommand{\tpmod}[1]{{\@displayfalse\pmod{#1}}}
\makeatletter
\@namedef{subjclassname@2020}{\textup{2020} Mathematics Subject Classification}
\makeatother

\begin{document}

\title{Super-positivity via Dirichlet theta functions}

\author{Arshay Sheth}
\address[Arshay Sheth]{School of Mathematics, Tata Institute of Fundamental Research, Homi Bhabha Road, Colaba, 
Mumbai - 400005, India}

\email{asheth@math.tifr.res.in}

\author{Wayne Voo}
\address[Wayne Voo]{Mathematics Institute, Zeeman Building, University of Warwick, Coventry, CV4 7AL, United Kingdom}
\email{wayne.voo@warwick.ac.uk}
\urladdr{}





\dedicatory{}

\begin{abstract}
Let $\chi$ be a quadratic Dirichlet character and let $\Lambda(\chi, s)$ be its associated completed $L$-function.  The super-positivity property, introduced by Yun--Zhang (2017), asserts that $\Lambda^{(m)} \left( \frac{1}{2}, \chi \right) \geq 0$ for all $m \geq 0$. We show that a certain positivity property of its  Dirichlet theta function $\theta(t, \chi)$, namely that $\theta(t, \chi)>0$ for all $t \geq 1$, implies the super-positivity property for $\chi$. We then investigate a refined form of the above property: we conjecture that the sequence $\{ \Lambda^{(2m)}\left(\frac{1}{2}, \chi \right )\}_{ m \geq 0}$ is always eventually increasing and has at most one turning point, and we prove these assertions again when $\theta(t, \chi)>0$ for $t \geq 1$. 
\end{abstract}

\maketitle

\section{Introduction}

Let $\chi$ be a primitive quadratic Dirichlet character and let $L(s, \chi)$ be its associated  $L$-function. The study of $L\left(\frac 1 2, \chi \right)$, and more generally the study of $L$-functions at the central point (as in, for example, the Birch--Swinnerton-Dyer conjecture),  plays an important role in number theory. 
Since $\chi$ is quadratic, $L\left(\frac 1 2, \chi \right) \in \R$ and so it is natural to investigate its sign. The Generalised Riemann Hypothesis (GRH) implies that $L\left(\frac 1 2, \chi \right)  \geq 0$. Indeed, the Euler product representation $$L(s,\chi)=\prod_{p} (1-\chi(p)p^{-s})^{-1} \text{ for } \Re(s)>1 $$ implies that $L(s,\chi)>0$ for $s \in \R_{>1}$,  so if $L\left(\frac 1 2, \chi \right)  <0$ then the Intermediate Value Theorem implies that $L(s, \chi)$ has a zero in $(\frac 1 2, 1+\epsilon]$ for some $\epsilon>0$, contradicting GRH. In particular, the weaker assumption that $L(s,\chi)$ has no real zeroes when $s > \frac 1 2$ also implies $L\left(\frac 1 2, \chi \right)  \geq 0$. 

Setting 
\[
\varepsilon (\chi) =
\begin{cases}
    0 \qquad \text{for } \chi \text{ even}, 
    \\
    1 \qquad \text{for } \chi \text{ odd}.
\end{cases}
\]
 and $\Gamma_{\mathbb{R}}(s)=\pi^{-\frac{s}{2}}\Gamma(\frac{s}{2}),$  the \textit{completed Dirichlet \(L\)-function} is defined by
\[
\Lambda(s,\chi):= N^{\frac{s}{2}}\Gamma_{\mathbb{R}}(s+\varepsilon(\chi))L(s,\chi),
\]
where $N$ is the conductor of $\chi$. Since $\Gamma(s)>0$ for $s \in \R_{> 0}$, GRH therefore implies that $\Lambda\left(\frac 1 2, \chi \right) \geq 0$.  In \cite{YZ17}, Yun--Zhang introduced the notion of super-positivity of $L$-functions, which extends this positivity property at the central point to all derivatives of $\Lambda(s, \chi)$.

\begin{definition}[{Yun--Zhang~\cite[Appendix B]{YZ17}}] \label{SP}
We say that $\Lambda(s, \chi)$ satisfies the super-positivity property if $\Lambda^{(m)} \left( \frac{1}{2}, \chi \right) \geq 0$ for all $m \geq 0$. 
\end{definition}

\begin{remark} \label{gen}
 Definition \ref{SP} can be extended to any $L$-function whose Dirichlet series coefficients are real numbers; these include, for example, $L$-functions attached to elliptic curves. More precisely,  if  $\Lambda(s, \pi)$ denotes the completed $L$-function attached to a self-dual cuspidal automorphic representation $\pi$ of $\mathrm{GL}_{n}(\mathbb{A}_{\Q})$, the super-positivity property in this case asserts that if $\Lambda^{(m)} \left( \frac{1}{2}, \pi \right) \geq 0$ for all $m \geq 0$. 
\end{remark}

Yun--Zhang show (\cite[Proposition B.1]{YZ17}) that GRH implies the super-positivity property for all $L$-functions in the setting of Remark \ref{gen}, so the super-positivity property is conjecturally true. 

\begin{remark}\label{odd even lambda}
 Since $\Lambda(\chi, s)$ satisfies the functional equation $\Lambda(\chi, s)=\Lambda(\chi, 1-s)$, it follows that $\Lambda^{ (m)} \left( \frac{1}{2}, \chi \right)= (-1)^m \Lambda^{ (m)} \left( \frac{1}{2}, \chi \right)$ for all $m \geq 0$. Thus, $\Lambda^{ (m)} \left( \frac{1}{2}, \chi \right)=0$ for $m$ odd and the super-positivity property is asserting that $\Lambda^{ (m)} \left( \frac{1}{2}, \chi \right) \geq 0$ for $m$ even (we in fact expect $\Lambda^{ (m)} \left( \frac{1}{2}, \chi \right)>0$ in these cases). More generally,   in the setting of Remark \ref{gen}, the super-positivity assertion is always automatically satisfied for one parity depending on the sign of its functional equation. 
 

\end{remark}

Writing the Taylor expansion of $\Lambda(s, \chi)$ around $\frac 1 2$  to be $\Lambda(s,\chi) = c_0+ (s-\frac{1}{2})c_1 + (s-\frac{1}{2})^2 c_2 + \dots$,  super-positivty asserts that $c_i \geq 0$ for all $i$. Since not all of the $c_i's$ are 0, super-positivity implies that $\Lambda(s,\chi)>0$ for $s \in \R_{>\frac{1}{2}}$.  We summarise the discussion so far via the following flowchart. 

\begin{figure}[htbp]
    \centering
    \begin{tikzpicture}[
        block/.style={
            draw,
            rectangle,
            minimum height=1.5em,
            inner sep=8pt,
            thick,
            align=center
        },
        implies/.style={
            double,
            double distance=2pt,
            -{Stealth[scale=1.2]},
            thick,
            shorten <=2mm,  
            shorten >=2mm   
        }
    ]

        \node[block] (grh) {GRH};
        \node[block, right=1.3cm of grh] (super) {Super-positivity};
        \node[block, right=1.3cm of super] (real) {no real zeros when $s>\frac 1 2$};
        \node[block, right=1.3cm of real] (siegel) {no Siegel zeros};
        \node[block, below=0.9cm of real] (lfunc) {$L\left(\frac{1}{2}, \chi\right) \geq 0$};

        \draw[implies] (grh) -- (super);
        \draw[implies] (super) -- (real);
        \draw[implies] (real) -- (siegel);
        \draw[implies] (real) -- (lfunc);

    \end{tikzpicture}
    \caption{Flowchart of logical implications related to the super-positivity property.}
\end{figure}
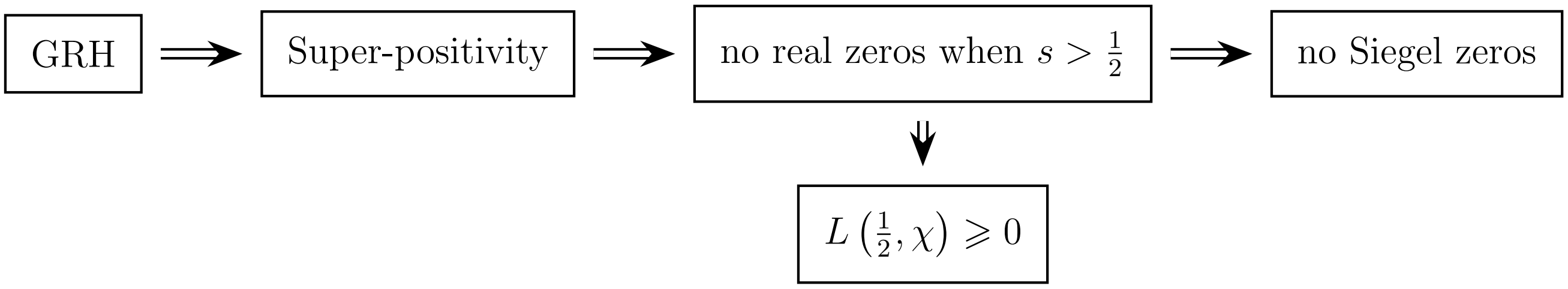

The goal of our paper is to study the super-positivity property for quadratic Dirichlet $L$-functions via the associated Dirichlet theta function $\theta(t, \chi)$, exploiting the connection between $\Lambda(s, \chi)$ and $\theta(t, \chi)$ via the Mellin transform. For $t>0$, we recall that  $\theta(t, \chi)$ is defined by

    \[ \theta(t, \chi)  = \begin{cases} 
          \sum_{m=1}^\infty \chi(m) \exp(-\pi t m^2 / N) & \text{ if } \chi \text{ is even } \\
         \sum_{m=1}^{\infty} m \chi(m) \exp{(- \pi t m^2 /N)} &  \text{ if } \chi \text{ is odd } 
       \end{cases}
    \]

 The Dirichlet theta function $\theta(t, \chi)$ satisfies the automorphy relation   

  \[ \theta(t, \chi)  = \begin{cases} 
          \frac{1}{\sqrt t} \cdot \theta(1/t, \chi) & \text{ if } \chi \text{ is even } \vspace{2mm} \\
        \frac{1}{t^{\frac 3 2}} \cdot \theta(1/t, \chi) &  \text{ if } \chi \text{ is odd } 
       \end{cases}
    \]

\begin{theorem} \label{mainthm1}
    Let $\chi$ be a primitive quadratic Dirichlet character. Then $\Lambda(s,\chi)$ satisfies super-positivity if \(\theta(t,\chi)>0\) for \(t\geq1\). 
\end{theorem}

We now make a few remarks about the assumption of $\theta(t,\chi)>0$ for $t \geq 1$  in Theorem \ref{mainthm1}, which has been studied, for example, by Louboutin--Munsch in \cite{LM13} and by Sarnak in \cite{Sar11}. It is conjectured that there are infinitely many Dirichlet theta functions which satisfy this property.  Numerical data (see \cite[page 648]{LM13}) reveals that the property holds for at least 6165194 of 30396324 possible positive fundamental discriminants up to $10^8$; Theorem \ref{mainthm1} applies in all these cases, and yields an unconditional verification of the super-positivity property. 
Note also that the automorphy relation for $\theta(t, \chi)$ implies that  $\theta(t,\chi)>0$ for $t \geq 1$ is equivalent to  $\theta(t,\chi)>0$ for $t > 0$.  
\\
We remark that a related property, namely the property of a quadratic character to have everywhere non-negative partial sums, has received much interest in recent literature; see, for example,
the work of Kalmynin \cite{Kal23} and the forthcoming work of Angelo--Soundararajan--Xu. In general, neither property implies the other. For example, the identity $\sum_{d \in (\mathbb Z/N \mathbb Z)^{\times}} \chi(d)=0$ implies that an even character never has everywhere non-negative partial sums, but \cite{LM13} give explicit examples of even characters whose theta function is everywhere positive.  However, Louboutin--Munsch show that if an even quadratic character has strictly positive partial sums on a specified initial range, then its associated theta function is everywhere positive (see \cite[Theorem 4]{LM13});  this fact yields a simple
way to check the assumption in Theorem \ref{mainthm1} and is indeed the basis of the numerical evidence
referenced above. 
\\
Finally, we remark that there are numerical methods to test super-positivity directly by showing that, for $s=\sigma+it$,  $\Lambda(s,\chi) \ne 0$ in the region where
$1/2 < \sigma < 1, \; |t| \leq \sigma - 1/2$ (see \cite[Theorem 1.2]{GH17}). The virtue of Theorem \ref{mainthm1} is that, by connecting super-positivity directly to the positivity properties of $\theta(t, \chi)$, we obtain an alternative  method/criterion to check  super-positivity (see also \cite[Remark B.5]{YZ17}).

\subsection{A refined conjectural property}

The super-positivity property can be regarded as an assertion that every number in the sequence $\{ \Lambda^{(2m)}\left(\frac{1}{2}, \chi \right )\}_{ m \geq 0}$ is non-negative. It is natural to investigate in further depth how this sequence of real numbers behaves. We present the following conjecture, which provides more detailed information about the sequence  $\{ \Lambda^{(2m)}\left(\frac{1}{2}, \chi \right )\}_{ m \geq 0}$. 

\begin{conjecture} \label{conj}
Let $\chi$ be a primitive quadratic Dirichlet character. Then the sequence 
\begin{enumerate}

\item $\{ \Lambda^{(2m)}\left(\frac{1}{2}, \chi \right )\}_{ m \geq 0}$ is eventually increasing.

\item  $\{ \Lambda^{(2m)}\left(\frac{1}{2}, \chi \right )\}_{ m \geq 0}$ has at most one turning point.

\end{enumerate}
\end{conjecture}

 As a first piece of evidence, we can prove our conjecture once again under the assumption that  \(\theta(t,\chi)>0\) for \(t\geq1\).

\begin{theorem} \label{mainthm2}
Let $\chi$ be a primitive quadratic Dirichlet character such that $\theta(t, \chi)>0$ for $t \geq 1$. Then the sequence 
\begin{enumerate}

\item  $\{ \Lambda^{(2m)}\left(\frac{1}{2}, \chi \right )\}_{ m \geq 0}$ is eventually increasing. 

\item $\{ \Lambda^{(2m)}\left(\frac{1}{2}, \chi \right )\}_{ m \geq 0}$ has at most one turning point. 

\end{enumerate}
\end{theorem}

We now discuss the numerical evidence for Conjecture \ref{conj}.  We verified the claims in the conjecture for the first ten derivatives for all quadratic Dirichlet $L$-functions of conductor up to 10,000. The numerical code that we used for our conjecture can be found at GitHub \href{https://github.com/arshaysheth24/Code-for-computing-derivatives-of-quadratic-Dirichlet-L-functions-at-the-central-point}{here}. We observed that there are examples of Dirichlet characters for which there is indeed a turning point, although we did not find any example once the conductor exceeds 97.  
\\
It would be an interesting project to investigate further properties of the sequence $\{ \Lambda^{(2m)}\left(\frac{1}{2}, \chi \right )\}_{ m \geq 0}$, such as its growth rate. For example, the proof of Theorem \ref{mainthm2} shows that, under the stated assumption $\theta(t, \chi)>0$ for $t \geq 1$, the sequence grows at least exponentially in $m$. It would be interesting to establish (even conditionally) an asymptotic growth of this sequence in general. 

\subsection*{AI Disclosure} ChatGPT was used to assist in developing the Python code for the numerical evidence for Conjecture \ref{conj}; the mathematical content in the paper is due to the authors.

\newpage 
\section{Proofs of Theorem \ref{mainthm1} and Theorem \ref{mainthm2} }

\begin{proof}[Proof of Theorem \ref{mainthm1}]

When $\chi$ is even, the integral representation of the completed Dirichlet $L$-function $\Lambda(s,\chi)$ (see, for example, \cite[Equation 2]{LM13}) yields

\begin{equation} \label{eqnzero}
\Lambda(s,\chi)=\int_1^\infty \theta(t,\chi) (t^{\frac{s}{2}}+t^{\frac{1-s}{2}}) \frac{dt}{t}
\end{equation}

It follows after taking the $m$-th derivative in \eqref{eqnzero} that
\begin{equation} \label{eqn1}
    \Lambda^{(m)}(s, \chi)= \int_1^\infty \theta(t, \chi) (\log \sqrt t)^m(t^{s/2}+(-1)^mt^{(1-s)/2}) \frac{dt}{t}. 
\end{equation}    

Substituting $s=\frac{1}{2}$ in \eqref{eqn1} yields

\begin{equation} \label{eqn2}
    \Lambda^{ (m)} \left( \frac{1}{2}, \chi \right) = \int_1^\infty \theta(t, \chi) (\log \sqrt t)^m(t^{-3/4}+ (-1)^m t^{-3/4}) dt.
\end{equation}

The assertion in Theorem \ref{mainthm1} now follows from \eqref{eqn2} and the hypothesis that $\theta(t, \chi)>0$ for $t \geq 1$. The case when $\chi$ is odd follows similarly, using the integral representation

\begin{equation} \label{eqnzero2}
N^{1/2} \Lambda(s,\chi)=\int_1^\infty \theta(t,\chi) (t^{-\frac{s}{2}}+t^{\frac{s-1}{2}}) dt. 
\end{equation}

\end{proof}

\begin{proof}[Proof of Theorem \ref{mainthm2}]
Again, we assume that $\chi$ is even since the proof below readily adapts to the odd case using \eqref{eqnzero2}. Define 
$$
f_\chi(x)= 2 \int_1^\infty \theta(t, \chi) \left(\frac{1}{2} \log t \right)^{2x} t^{-3/4}  dt,
$$
so \eqref{eqn2} implies that  $f_\chi(m)= \Lambda^{ (2m)} \left( \frac{1}{2}, \chi \right)$ for $m \geq 0$. To show (1) and (2), it suffices to show that $f_\chi(x) \to \infty$ as $x \to \infty$ and $f_\chi^{''}(x)>0$ for all $x>0$. Indeed, an application of the mean value theorem shows that any real valued function which on $(0, \infty)$ which satisfies these two properties is eventually increasing; furthermore, any function whose second derivative is positive on $(0, \infty)$ has at most one turning point. 

Since the integrand is nonnegative over the integration range by hypothesis, one deduces
\begin{equation}\label{eventually infinity}
\int_1^\infty \theta(t, \chi) \left(\frac{1}{2} \log t \right)^{2x} t^{-3/4}  dt \geq \int_{e^3}^{e^4}\theta(t, \chi) \left(\frac{1}{2} \log t \right)^{2x} t^{-3/4} dt. 
\end{equation}

By continuity of $\theta(t,\chi)$ and $t^{-3/4}$ on $[e^3,e^4]$ and again using the assumption $\theta(t, \chi) \geq 0$ for $t \geq 1$, there exists $M\in \R_{>0}$ such that $\theta(t,\chi)t^{-3/4} \geq M$ on this interval. Hence, from \eqref{eventually infinity}, it follows that 
$$
\int_1^\infty \theta(t, \chi_D) \left(\frac{1}{2} \log t \right)^{2x} t^{-3/4}  dt \geq (e^4-e^3)M \left( \frac 3 2 \right)^{2x},
$$

which tends to infinity as $x \rightarrow \infty$. Therefore, $f_\chi(x) \to \infty$ as $x \to \infty$. Finally, taking the second derivative yields 
$$
f_\chi(x)''= 2 \int_1^\infty \theta(t, \chi) (2 \log \log \sqrt t)^2 \log (\sqrt t)^{2x}   t^{-3/4}  dt >0.  \qedhere  
$$

\end{proof}

\appendix

\section{Super-positivity for the Riemann zeta function}
It was mentioned on \cite[page 906]{YZ17} that super-positivity for the Riemann zeta function was proven by P\'olya with a reference to the article \cite{CNV86} (which in-turn references \cite{Pol27}). Such a claim was also made on \cite[page 47]{GH17}. However, we could not find a precise proof of this assertion in any of the references, so we give a short self-contained proof below. The proof is similar to the proof of Theorem \ref{mainthm1} (see also \cite[Theorem 18]{KS26}), relying on an integral representation with positive kernel. We let $\Lambda(s)=\frac{1}{2}(s)(s-1)\pi^\frac{-s}{2} \Gamma(\frac{s}{2}) \zeta(s)$ denote the completed Riemann zeta function. 

\begin{theorem}
We have that $\Lambda^{(m)} \left( \frac{1}{2} \right) \geq 0$ for all $m \geq 0$. 
\end{theorem}

\begin{proof}
As explained in Remark \ref{odd even lambda}, $\Lambda^{(m)} \left( \frac{1}{2} \right) = 0$ for $m$ odd.
We will  use the following integral representation of the completed Riemann Zeta Function (see, for example,  \cite[Equation (10.1.3)]{Tit86} ): namely,  $\Xi(t):= \Lambda \left (\frac 1 2 + it \right) = 2\int_0^\infty \cos(tx) \Phi(x)dx $, where $ \Phi(x)=\sum_{n=1}^\infty \left(4 \pi^2 n^4 e^{\frac{9x}{2}}- 6 \pi n^2 e^{\frac{5x}{2}}\right)e^{-\pi n^2 e^{2x}}$. Note that $\Phi(x) \geq 0$ for all $x \geq 0$. Taking the $m$-th derivative yields $\Xi^{(m)}(t)= 2\int_0 ^\infty f_m(x) \Phi(x) dx$, where 
    \[f_m(x)=
    \begin{cases}
        (-1)^{(m+1)/2} x^m \sin(tx),  &\text{for }m \text{ odd,}
        \\
        (-1)^{m/2} x^m \cos (tx), &\text{for } m \text{ even.}
    \end{cases}
    \]

On the other hand, we have $
    \Xi^{(m)}(t) = (i)^m \Lambda^{(m)}(\frac{1}{2} + it) = (-1)^{m/2} \Lambda^{(m)}(\frac{1}{2} + it)$.  Hence, substituting $t=0$ for $m$ even, $\Lambda^{(m)}(\frac{1}{2})=  (-1)^{m/2} \Xi^{(m)}(0)= 2\int_0^\infty x^m \Phi(x) dx \geq 0$.  
\end{proof}


\begin{thebibliography}{BCDTO1}

\bibitem[CNV86]{CNV86}
George Csordas, Timothy S. Norfolk and Richard S. Varga, 
\href{https://www.jstor.org/stable/2000378}{\emph{The {R}iemann hypothesis and the {T}ur\'an inequalities}}.
 Trans. Amer. Math. Soc. \textbf{296} (1986), no.2, 521--541, 1986.

\bibitem[GH17]{GH17}
Dorian Goldfeld and Bingrong Huang, \href{https://www.impan.pl/en/publishing-house/banach-center-publications/all/118/0/113098/super-positivity-of-a-family-of-l-functions}{\emph{Super-positivity of a family of L-functions}}, Banach Center Publications \textbf{118} (2019), 45–93.


\bibitem[LM13]{LM13}
Stéphane R. Louboutin and Marc Munsch, \href{https://publi.math.unideb.hu/paper/1837/download/10_5486_PMD_2013_5670.pdf}{\emph{On positive real zeros of theta and \(L\)-functions associated with real, even and primitive characters}}, Publicationes Mathematicae, Vol. 83, No. 4, 01.12.2013, p. 643-665.

\bibitem[Kal23]{Kal23}
Alexander Kalmynin, \href{https://londmathsoc.onlinelibrary.wiley.com/doi/10.1112/mtk.12179}{\emph{Quadratic characters with positive partial sums}}, 
Mathematika \textbf{69} (2023), no. 1, 90–99.



\bibitem[KS26]{KS26}
Shin-ya Koyama and Arshay Sheth, \href{https://comptes-rendus.academie-sciences.fr/mathematique/articles/10.5802/crmath.846/}{\emph{Chebyshev's bias for modular forms}}, Comptes Rendus. Mathématique, Volume 364 (2026),  pp. 589-598. 


\bibitem[Pol27]{Pol27}
George P\'olya, 
\emph{\"Uber die algebraisch-funktionentheoretische Untersuchungen von J.L.W.V. Jensen}, 
Kgl. Danske Md. Sel. Math. Fys. Medd. \textbf{7} (1927), 3-33.

\bibitem[Sar11]{Sar11}
Peter Sarnak, 
\newblock \href{http://publications.ias.edu/sarnak/paper/511, 2011.}{\emph{Letter to Eitan Bachmat on positive definite L-Functions}}, 2011. 

\bibitem[Tit86]{Tit86}
Edward Charles Titchmarsh, {\emph{The theory of the Riemann zeta function}}, second edition, 
The Clarendon Press, Oxford University Press, New York, 1986. 

\bibitem[YZ17]{YZ17}
Zhiwei Yun and Wei Zhang, \href{https://annals.math.princeton.edu/2017/186-3/p02}{\emph{Shtukas and the Taylor expansion of $L$-functions}}, Ann. of Math. (2) \textbf{186} (2017), no.3, 767-911. 

\end{thebibliography}
\end{document}